\documentclass[10pt]{smfart}

\RequirePackage[T1]{fontenc}
\RequirePackage{amsfonts,latexsym,amssymb}
\RequirePackage[frenchb]{babel}
\addto\extrasfrenchb{\bbl@nonfrenchitemize
\bbl@nonfrenchspacing}

\RequirePackage{mathrsfs}
\let\mathcal\mathscr

\theoremstyle{plain}
\numberwithin{equation}{section}
\newtheorem{prop}[equation]{\propname}
\newtheorem{theo}[equation]{\theoname}

\theoremstyle{definition}
\theoremstyle{remark}

\usepackage[matrix,arrow]{xy}
\usepackage{url}
\usepackage{fge}
\usepackage{bbm}
\newcommand{\moins}{\mathbin{\fgebackslash}}
\usepackage[pagebackref]{hyperref}

\let\cal\mathcal

\def\paskunas{{Pa\v{s}k\={u}nas}}
\def\bPi{\boldsymbol\Pi}

\def\Q{{\bf Q}} \def\Z{{\bf Z}}
 
\def\C{{\bf C}}

\def\dual{{\boldsymbol *}}

\def\Qbar{\overline{\bf Q}}
\def\epsilon{\varepsilon}

\def\ainf{{\bf A}_{{\rm inf}}}
\def\bst{{\bf B}_{{\rm st}}}
\def\bcris{{\bf B}_{{\rm cris}}}

\def\Bcris{{\mathbb B}_{{\rm cris}}}

\def\piqp{{\bf P}^1}
\def\wotimes{\,\widehat\otimes\,}
\def\GL{{\bf GL}}

\def\matrice#1#2#3#4{{\big(\begin{smallmatrix}#1&#2\\ #3&#4\end{smallmatrix}\big)}}

\newcommand{\eet}{\operatorname{\acute{e}t} }
\newcommand{\proet}{\operatorname{pro\acute{e}t} }

\def\lra{\leftrightarrow}

\def\O{{\cal O}}

\begin{document}

\title[$p$-adic local Langlands correspondence]{On the geometrization of the $p$-adic local Langlands correspondence}

\author{Pierre Colmez}
\address{CNRS, IMJ-PRG, Sorbonne Universit\'e, 4 place Jussieu, 75005 Paris, France}
\email{pierre.colmez@imj-prg.fr}

\author{Gabriel Dospinescu}
\address{CNRS, LMBP, Universit\'e Clermont-Auvergne, Campus des C\'ezeaux, 3 place Vasarely,
63170 Aubi\`ere, France}
\email{gabriel.dospinescu@uca.fr}

\author{Wies{\l}awa Nizio{\l}}
\address{CNRS, IMJ-PRG, Sorbonne Universit\'e, 4 place Jussieu, 75005 Paris, France}
\email{wieslawa.niziol@imj-prg.fr}

\begin{abstract}
We survey results related to our geometrization of a part of the $p$-adic local Langlands
correspondence for ${\bf GL}_2(\Q_p)$.
\end{abstract}

\maketitle

\section*{Introduction}
The $p$-adic local Langlands correspondence for ${\bf GL}_2(\Q_p)$ has a
very satisfactory shape but there is no other group for which the situation is as well 
understood (not even ${\bf GL}_3(\Q_p)$ or ${\bf GL}_2(F)$ for $[F:\Q_p]<\infty$).
There exist though several competing approaches for remedying this situation. Here we will review 
our proposal from~\cite{CDN1} 
and some related research.
\section{The $p$-adic local Langlands correspondence for ${\bf GL}_2(\Q_p)$}
\subsection{Notation}
In what follows we fix a finite extension $L$ of $\Q_p$ (this will be the field
of coefficients for our representations of groups), with ring of integers
$\O_L$ and residue field $k_L$.

Let ${\rm Gal}_{\Q_p}\supset{\rm W}_{\Q_p}$ denote
the absolute Galois group of $\Q_p$ and the Weil group of $\Q_p$.
Local class field theory gives an identification between continuous
characters of ${\rm W}_{\Q_p}$ and continuous characters of $\Q_p^\dual$,
and a character of ${\rm W}_{\Q_p}$ extends to ${\rm Gal}_{\Q_p}$ if
and only if it is unitary (i.e., takes values in $\O_L^\dual$);
under this identification the cyclotomic character $\chi:{\rm Gal}_{\Q_p}\to \O_L^\dual$
corresponds to $x\mapsto x|x|$.

We set $G:={\bf GL}_2(\Q_p)$.
\subsection{The correspondence}
Let 
${\rm Rep}_L\,G$ be the category of continuous representations
of $G$ on $L$-Banach spaces $\Pi$, that are unitary (there is a $G$-invariant 
norm defining the topology
on $\Pi$), of finite length as topological $G$-modules (this implies that $\Pi$ is
admissible, i.e.~if $\Pi_0$ is the unit ball for a $G$-invariant norm, the action
of $G$ on $k_L\otimes_{\O_L}\Pi_0$ is locally constant and, if $K$ is an open subgroup of
$G$, the $K$-fixed vectors are finite dimensional over $k_L$).

Let 
${\rm Rep}_L\,{\rm Gal}_{\Q_p}$ be the category of finite dimensional $L$-modules
with a continuous and linear action of ${\rm Gal}_{\Q_p}$. If
$n\geq 1$, let ${\rm Rep}^{(n)}_L\,{\rm Gal}_{\Q_p}$ be the subset of representations
of dimension~$n$.

\begin{theo}\phantomsection\label{cdn1} There exists:

{\rm (i)} 
a functor $\Pi\mapsto{\bf V}(\Pi)$, 
from ${\rm Rep}_L\,G$ to $ {\rm Rep}_L\,{\rm Gal}_{\Q_p}$,

{\rm (ii)} 
a correspondence $V\mapsto\bPi(V)$
from ${\rm Rep}^{(2)}_L\,{\rm Gal}_{\Q_p}$ to $ {\rm Rep}_L\,G$ with:

\quad $\bullet$ ${\bf V}(\bPi(V))=V$,

\quad $\bullet$ $\bPi(V)$ has no finite subobject,

\quad $\bullet$ $\bPi(V)$ is maximal for these two properties.
\end{theo} 
The construction of the functor ${\bf V}$ goes through the theory of $(\varphi,\Gamma)$-modules
introduced~\cite{Fo91} by Fontaine to classify $L$-representations of ${\rm Gal}_{\Q_p}$,
and is rather straightforward.
What is delicate is to show that one gets finite dimensional
representations of ${\rm Gal}_{\Q_p}$ and that any representation of dimension $2$
is in the image.

\subsection{Properties of the correspondence}
The above correspondence enjoys a number of additional properties:

\vskip1mm
$\bullet$ 
Compatibility with class field theory: the central character of $\bPi(V)$
and the determinant of $V$ are related by the formula:
$\omega_{\bPi(V)}=\chi^{-1}\det V.$ 

\vskip1mm
$\bullet$ $V\mapsto\bPi(V)$ induces a bijection between $2$-dimensional,
absolutely irreductible $V$'s, and
supercuspidal $\Pi$'s (i.e., absolutely irreducible, and not a Jordan-H\"older
factor of the induction of a character of the Borel).

If $V$ is not irreducible, neither is $\bPi(V)$; more precisely:
$$V^{\rm ss}=\delta_1\oplus\delta_2 \Rightarrow
\bPi(V)^{\rm ss}\sim {\rm Ind}_B^G(\delta_2\otimes\delta_1\chi^{-1})^{\rm ss}
\oplus {\rm Ind}_B^G(\delta_1\otimes\delta_2\chi^{-1})^{\rm ss}$$

Moreover, all irreducible $\Pi$'s appear as a Jordan-H\"older factor of a  $\bPi(V)$.

\vskip1mm
$\bullet$ 
Relation to the classical local Langlands correspondence:
$\bPi(V)^{\rm alg}\neq 0$ if and only if  $V$ is de Rham, 
with distinct Hodge-Tate weights $a<b$, and we have:
$$\bPi(V)^{\rm alg}={\rm LL}(D_{\rm pst}(V))\otimes {\rm Sym}^{b-a-1}\otimes{\det}^a$$
Moreover, there is also a recipe to recover ${\rm LL}(D_{\rm pst}(V))$ from $\bPi(V)$ if $a=b$.

\vskip1mm
$\bullet$ Compatibility with a mod~$p$ local Langlands correspondence: in a semi-simplified
version, this mod~$p$ correspondence gives a bijection ${\cal B}\lra \overline\rho_{\cal B}$
between the blocks of ${\rm Rep}_{k_L}\,G$ and the 
$2$-dimensional, semi-simple, $k_L$-representations of ${\rm Gal}_{\Q_p}$. 

\vskip1mm
$\bullet$ Compatibility with global correspondence:
$\bPi(V)$ appears in the completed cohomology of towers of modular curves.

\vskip2mm
One useful application of the correspondence is the following result explaining
how to recover $\Pi$ from its locally analytic vectors. 
\begin{prop}\phantomsection\label{cd}
 {\rm(\cite{CD})}
The dual $\Pi^\dual$
of $\Pi$ is the subspace of $(\Pi^{\rm an})^{\dual}$ of vectors which are
bounded under the action of $G$.
\end{prop}
A natural question (Emerton) is whether this relation is true for 
finite length
unitary admissible representations of ${\bf GL}_n(F)$, $[F:\Q_p]<\infty$
(or more generally any reductive group in place of ${\bf GL}_n$).

 \subsection{Some history}
The above results are the combination of a lot of research of different people.
We give a brief summary of the main steps.

\vskip1mm
$\bullet$ Barthel-Livn\'e (1994-95) partially classified~\cite{BL94} $k_L$-representations 
of ${\bf GL}_2(F)$; this was 
completed~\cite{Br1} by Breuil (2003) for ${\bf GL}_2(\Q_p)$. \paskunas\ computed~\cite{Pas1}
 the blocks for ${\bf GL}_2(\Q_p)$. 

\vskip1mm
$\bullet$ Schneider-Teitelbaum (2000-03) developped~\cite{ST1,ST2,ST3,ST4}
the foundations of the theory of $p$-adic representations of $p$-adic groups
(admissibility, locally analytic and algebraic vectors).
The missing Schur Lemma was provided~\cite{DS} by
Dospinescu-Schraen (2012) building on a result of Ardakov-Wadsley~\cite{AW} .

\vskip1mm
$\bullet$ Breuil (2000 -- ...)
understood~\cite{Br4}  what could be a $p$-adic Local Langlands correspondence 
for de Rham representations: a
de Rham representation of ${\rm Gal}_{\Q_p}$ is classified by the data of
a filtered representation of the Weil-Deligne group of $\Q_p$; the classical
local Langlands correspondence encodes the Weil-Deligne part but says nothing
about the filtration, and Breuil had the insight that possible filtrations 
on a Weil-Deligne representation should correspond
to completions of a locally algebraic representation built from the classical one
and an algebraic representation encoding the jumps of the filtration.
He also made
precise conjectures~\cite{Br2,Br3} for crystalline and semi-stable $2$-dimensional 
representations of ${\rm Gal}_{\Q_p}$. 

\vskip1mm
$\bullet$ Emerton (2003 -- ...) 
introduced~\cite{Em05,Em06}
 the concept of completed cohomology (of towers of modular curves, Shimura varieties or
symmetric spaces), and formulated~\cite{Em06b} a conjectural local-global compatibility
for the tower of modular curves. 

\vskip1mm
$\bullet$ 
The above mentionned conjectures of Breuil were proven~\cite{principal,BB}
 by Colmez-Berger-Breuil (2004-05)
through the connection to the theory of $(\varphi,\Gamma)$-modules,
first in the semi-stable case (Colmez), followed immediately by the crystalline case (Berger-Breuil).
These results were extended to trianguline representations by Colmez.

\vskip1mm
$\bullet$ Colmez (2005 -- 2008) defined~\cite{gl2} the functor ${\bf V}$ and gave a direct construction
of $\bPi$ (Kisin had suggested a construction by deformation, using the density of crystalline
representations), allowing the study of locally analytic
and algebraic vectors. The surjectivity of $\bPi$ relies on the density of crystalline
representations which was not known at the time in some cases for $p=2$.
Dospinescu (2012) proved~\cite{Dosp1} the existence of an infinitesimal character and that this
infinitesimal character encodes the Hodge-Tate weights of the associated representation
of ${\rm Gal}_{\Q_p}$. The description of the locally algebraic vectors found
applications to the Fontaine-Mazur conjecture (Kisin~\cite{Ki-FM}, 
Emerton~\cite{Em08}, Pan~\cite{pan1}).

\vskip1mm
$\bullet$ Emerton (2008) proved in~\cite{Em08} big parts of his local-global compatibility conjecture;
this was later revisited and complemented by  
Caraiani-Emerton-Gee-Geraghty-\paskunas-Shin~\cite{6auteurs1}, 
Pan~\cite{pan1,pan2},
Colmez-Wang~\cite{adele, CW}. 

\vskip1mm
$\bullet$ \paskunas\ (2012) proved in~\cite{Pas} (with some restrictions for $p=2,3$)
that ${\bf V}$ applied to irreducible objects produces
representations of dimension~$2$ except for Jordan-H\"older factors of representations induced
from a character of the Borel where the dimension is~$1$ (or~$0$ if it is finite dimensional,
i.e., is a character).  He also described the projective envelope of 
(the dual) of an irreducible $k_L$-representation $\pi$ of ${\bf GL}_2(\Q_p)$ in Galois terms
(i.e., in terms of the universal deformation of the representation of ${\rm Gal}_{\Q_p}$
associated to the block of $\pi$; in retrospect this amounts to a local categorification
of the correspondence).

\vskip1mm
$\bullet$ Colmez-Dospinescu-\paskunas\ (2014) treated in~\cite{CDP} the remaining cases (for $p=2,3$)
and proved compatibility with local class field theory. See also~\cite{PT} for extra cases
of local categorification.

\vskip1mm
$\bullet$ Colmez gave a recipe~\cite{poids} to extract, from
the unitary representation of ${\bf GL}_2(\Q_p)$, the classical representation
attached to a de Rham representation with equal Hodge-Tate weights.

\vskip1mm
$\bullet$ Dotto-Emerton-Gee recently announced a full categorification of the correspondence
(i.e., an equivalence of categories between residually finitely presented representations
of ${\bf GL}_2(\Q_p)$ and coherent sheaves on the Emerton-Gee stack parametrizing
$(\varphi,\Gamma)$-modules of rank~$2$).

 \subsection{Other groups}
As one can see from the above mentionned results, 
the case of ${\bf GL}_2(\Q_p)$ is remarkably well understood.  There is
no group not directly connected to ${\bf GL}_2(\Q_p)$ 
(such as ${\bf SL}_2(\Q_p)$ or ${\bf PGL}_2(\Q_p)$)
for which we know what to expect.
There have been different approaches to attack the problem.

\vskip1mm
$\bullet$ Caraiani-Emerton-Gee-Geraghty-\paskunas-Shin (2015)
construct~\cite{6auteurs} many $\bPi(V)$'s for many $V$'s 
(with $V\in{\rm Rep}_L^{(n)}{\rm Gal}_F$, $[F:\Q_p]<\infty$,
$\bPi(V)\in{\rm Rep}_L{\bf GL}_n(F)$), using global methods (completed cohomology of towers
of Shimura varieties), but $\bPi(V)$ could depend on global choices.

\vskip1mm
$\bullet$ Breuil and collaborators (Ding, Hellmann, Herzig, Hu, Morra, Schraen, Wang) 
construct~\cite{Br5,BHg,BD1,BD2,BHS}
 parts of what should
be $\bPi(V)^{\rm an}$, hopefully big enough to encode $V$ 
(crystalline or semi-stable case). A reality check is to verify that the
constructed representation appears in completed cohomology; this has applications
to the Fontaine-Mazur conjecture~\cite{BHS}.

\vskip1mm
$\bullet$ 
Emerton-Gee-Hellmann (2022) formulated in~\cite{EGH}
 a general categorification conjecture.

\vskip1mm
$\bullet$ Our approach~\cite{CDN1,CDN3,CDN4,CDN5} is to geometrize the existing correspondence 
for ${\bf GL}_2(\Q_p)$ by showing that it appears
in the $p$-adic \'etale 
cohomology of the Drinfeld tower (the $\ell$-adic \'etale
cohomology is known to encode the classical correspondence).
This gives a potential candidate for ${\bf GL}_n(F)$, with $[F:\Q_p]<\infty$.

\section{Geometrization of the correspondence}
 \subsection{The Drinfeld tower}
The Drinfeld upper half-plane $\Omega_{\rm Dr}:=\piqp\moins\piqp(\Q_p)$
is an analytic curve defined over $\Q_p$. It comes with an
 action of $G$ given by $\matrice{a}{b}{c}{d}=\frac{ax+b}{cx+d}$.

The Drinfeld
Tower~\cite{Drinfeld,Drinfeld1}
 $\Omega_{\rm Dr}\leftarrow{\cal M}_0\leftarrow{\cal M}_1\leftarrow\cdots {\cal M}_\infty$, 
is a tower of \'etale
$G$-equivariant coverings (cofinal in such coverings by Scholze-Weinstein)
 defined over $\breve\Q_p:=\widehat{\Q_p^{\rm nr}}$,
with an action of ${\rm W}_{\Q_p}$ (not ${\rm Gal}_{\Q_p}$). 
We have
$${\cal M}_0=\Omega_{\rm Dr}\times\Z,\quad \pi_0({\cal M}_{\infty,\C_p})=\Q_p^\dual$$
and, if $\check G=D^\dual$, where $D/\Q_p$ is the quaternion algebra, with ring
of integers $\O_D$, and uniformizer $\pi_D$,
$${\rm Aut}({\cal M}_\infty/\Omega_{\rm Dr})=\check G,\quad 
{\rm Aut}({\cal M}_n/{\cal M}_0)=\O_D^\dual/(1+\pi_D^n\O_D)$$

If $H^\dual$ is any reasonnable cohomology theory, we have an
action of $G\times \check G\times{\rm W}_{\Q_p}$ on 
$H^*({\cal M}_{\infty,\C_p}):=\varinjlim_n H^*({\cal M}_{n,\C_p})$.
For example,
the following result is classical.

\begin{theo}\phantomsection\label{cdn2}
If $\ell\neq p$,
$$\Qbar_\ell\otimes H^1_{\proet}({\cal M}_{\infty,\C_p},\Q_\ell(1))=
\bigoplus_M {\rm cosoc}(M)\otimes {\rm LL}(M)^\dual
\otimes {\rm JL}(M)$$
where the
sum\footnote{Actually, one has to quotient 
${\cal M}_{\infty,\C_p}$ by $\matrice{p}{0}{0}{p}\in\check G$ and restrict to $M$'s
with determinant $p$ on Frobenius for the statement to be correct; there is a way to
recover $H^1_{\proet}({\cal M}_{\infty,\C_p},\Q_\ell(1))$ from the cohomology of this quotient.} 
 is over $2$-dimensional, indecomposable, $\Qbar_\ell$-representations $M$ of ${\rm WD}_{\Q_p}$, 
${\rm cosoc}(M)$ is the irreducible quotient {\rm ($=M$ in general)}.
\end{theo}
\noindent{\bf Remark}\hskip2mm
Drinfeld did not restrict himself to ${\bf GL}_2(\Q_p)$ and constructed a tower corresponding
to ${\bf GL}_n(F)$ for any $n\geq 2$ and $[F:\Q_p]<\infty$.  The bottom of the tower is obtained
by removing from ${\bf P}^{n-1}$ all hyperplanes defined over $F$.
The analog of the above result is valid in this generality by works of
Drinfeld, Deligne, Carayol,  Harris-Taylor, Faltings,
Fargues, Dat, Boyer.

\smallskip
In summary the $\ell$-adic \'etale cohomology of the Drinfeld tower, for $\ell\neq p$,
 encodes simultaneously the classical local Langlands and Jacquet-Langlands correspondences.
An obvious idea to geometrize the $p$-adic local Langlands correspondence is to set $\ell=p$,
but this has scared people for a while because $p$-adic (pro-)\'etale cohomology of
$p$-adic analytic spaces is much more complicated than their $\ell$-adic one, for $\ell\neq p$.
For example, we have the following result (special case of Th.\,\ref{stein}):
\begin{prop}\phantomsection\label{cdn3}
If $Y$ is the open unit ball over $\Q_p$, then
$$H^1_{\proet}(Y_{\C_p},\Q_\ell(1))=\begin{cases}
0&{\text{if $\ell\neq p$,}}\\ 
\Omega^1(Y_{\C_p}) &{\text{if $\ell= p$.}}
\end{cases}$$
In particular, $H^1_{\proet}(Y_{\C_p},\Q_p(1))$ contains any $V\in{\rm Rep}_L\,{\rm Gal}_{\Q_p}$
with at least one Hodge-Tate weight $0$.
\end{prop}

 \subsection{Multiplicities}
Set $$H^1_{{\eet},\C_p}:=\varinjlim\nolimits_n H^1_{\eet}({\cal M}_{n,\C_p},L(1))$$
Say that {\it $V\in{\rm Rep}_L\,{\rm Gal}_{\Q_p}$ is OK} 
if it is $2$-dimensional, de Rham with Hodge-Tate weights $0$ and $1$, 
and the Weil-Deligne representation $D_{\rm pst}(V)$ is irreducible\footnote{A typical
example of OK representation is an irreducible representation appearing in the $p$-adic
\'etale cohomology
of the modular curve $X_0(N)$ with $p^2\mid N$ and not 
in that of $X_0(N)$ if $p^2\nmid N$.}.
Then the representation
${\rm LL}(V):={\rm LL}(D_{\rm pst}(V))$ of $G$ associated
by the classical local Langlands correspondence is supercuspidal, and
the representation 
${\rm JL}(V):={\rm JL}({\rm LL}(V))$ of $\check G$ associated
by the classical Jacquet-Langlands correspondence
is finite dimensional and irreducible.

\begin{theo}\phantomsection\label{cdn4}
{\rm(\cite{CDN1})}
Let $V$ be irreducible, of dimension~$\geq 2$.
 $${\rm Hom}_{W_{\Q_p}}(V,H^1_{{\eet},\C_p})=\begin{cases}
\bPi(V)^\dual\otimes {\rm JL}(V) &{\text{if $V$ is OK,}}\\ 0&{\text{otherwise.}}\end{cases}$$
\end{theo}

(In summary: 
$H^1_{{\eet},\C_p}$ contains the Galois representations that one would like with the multiplicity
that one would wish for.)

\vskip3mm
Say that {\it $\Pi\in {\rm Rep}_L\,G$ is OK} 
if $\Pi=\bPi(V)$ with $V$ OK, or if $\Pi$ is a twist of the Steinberg
representation by a locally constant
character.

\begin{theo}\phantomsection\label{cdn5}
{\rm (\cite{CDN5})}
Let $\Pi$ be irreducible.
 $${\rm Hom}_{G}(\Pi^\dual,H^1_{{\eet},\C_p})=\begin{cases}
{\bf V}(\Pi)\otimes {\rm JL}(\Pi^{\rm sm}) &{\text{if $\Pi$ is OK,}}\\ 0&{\text{otherwise.}}\end{cases}$$
\end{theo}

Both results have been generalized by Vanhaecke~\cite{VanH} to general Hodge-Tate weights (using non
constant local systems).

 \subsection{Families}
$H^1_{{\eet},\C_p}$ seems to be too big to have a reasonnable description. Set
 $$H^1_{{\eet},\Qbar_p}:=\varinjlim\nolimits_n 
L\otimes\big(\varprojlim_k\big(\varinjlim_{[K:\Q_p]<\infty}
 H^1_{\eet}({\cal M}_{n,K},(\O_L/p^k)(1))\big)\big)
\subset H^1_{{\eet},\C_p}$$

\begin{theo}\phantomsection\label{cdn6}
{\rm (\cite{CDN5})} We have a factorization
 $$H^1_{{\eet},\Qbar_p}\cong\bigoplus_M\Big(\bigoplus_{\cal B}\bPi(\rho_{M,{\cal B}})^\dual
\hskip-2mm\underset{R_{{\cal B},M}}{\otimes}\hskip-2mm
\rho_{M,{\cal B}}\Big)\underset{L}{\otimes}{\rm JL}(M),$$
where
$M$ runs through 
``irreducible $2$-dimensional representations of $W_{\Q_p}$'', and ${\cal B}$ 
through blocks of ${\rm Rep}_{k_L}\,G$,
$\rho_{M,{\cal B}}:{\rm Gal}_{\Q_p}\to {\bf GL}_2(R_{{\cal B},M})$ 
is the universal OK deformation of $\overline\rho_{\cal B}$ with
$D_{\rm pst}=M$, and
$\bPi(\rho_{M,{\cal B}})$ is  the representation of $G$ attached to $\rho_{M,{\cal B}}$ by 
the family version of the $p$-adic local Langlands correspondence.
\end{theo}

\noindent{\bf Remark}\hskip2mm
Let $M_{\rm dR}:=(\Qbar_p\otimes_{\Q_p^{\rm nr}}M)^{{\rm Gal}_{\Q_p}}$, 
a $2$-dimensional $L$-module. If
${\cal L}\subset M_{\rm dR}$ is a line, then, thanks to~\cite{CF},
 $$V_{M,{\cal L}}:={\rm Ker}((\bcris^+\otimes M)^{\varphi=p}\to\C_p\otimes
(M_{\rm dR}/{\cal L}))$$
is an OK representation of ${\rm Gal}_{\Q_p}$ and
all OK $V$'s are a $V_{M,{\cal L}}$ for a unique couple $(M,{\cal L})$.
One infers from this that $R_{{\cal B},M}$ is the ring of bounded analytic functions
on a nice open subset of $\piqp(M_{\rm dR})$ which implies that it is a product of PID's.

 \subsection{Proofs}
The starting point is that one can recover $H^1_{\eet}$ as the space of 
$v\in H^1_{\proet}$ whose orbit under $G$ is bounded.  
Hence we first compute $H^1_{\proet}$
which is easier because it involves rational $p$-Hodge theory instead of integral
$p$-Hodge theory.

\vskip2mm
${\cal M}_\infty$ is a (projective system of) curve and even a Stein curve (strictly increasing union
of affinoids). So, let
$Y/\C_p$ be a curve, and set
$\O:=\O(Y)$, $\Omega^1:=\Omega^1(Y)$, $H^1_{\proet}:=H^1_{\proet}(Y,\Q_p(1))$,
$H^1_{?}:=H^1_{?}(Y)$ if $?={\rm HK},{\rm dR}$.
In particular, $H^1_{\rm HK}$ is a $\Q_p^{\rm nr}$-module with an action of a Frobenius $\varphi$,
a monodromy operator $N$, and a 
\linebreak
(pro-)smooth action of ${\rm Gal}_K$ if $Y$ is defined
over $K$, and we have a Hyodo-Kato isomorphism
$\iota_{\rm HK}:\C_p\wotimes_{\Q_p^{\rm nr}}H^1_{\rm HK}(Y)\overset{\sim}\to H^1_{\rm dR}(Y)$.

As usual, a module with compatible actions of $\varphi$, $N$ and ${\rm Gal}_K$ can be turned into
a representation of ${\rm WD}_K$ and we have:
\begin{theo}\phantomsection\label{cdn7}
 {\rm(\cite{CDN2})}
If $\ell\neq p$ the representation of ${\rm WD}_{K}$ on 
$ H^1_{\proet}(Y,\Q_\ell(1))$ is "isomorphic"\footnote{With quotation
marks because they are not over the same field.} to that on $H^1_{\rm HK}(Y)$.
\end{theo}

\begin{theo}\phantomsection\label{stein}
 {\rm (\cite{CDN1,CDN2,CDN3})}
If $Y/\C_p$ is Stein, we have a functorial diagram of Fr\'echet spaces\label{cdn8}
$$\xymatrix@R=.5cm@C=.5cm{
0\ar[r] & \O/\C_p\ar[r]\ar@{=}[d] & H^1_{\proet}\ar[d] \ar[r] &
(\bst^+\wotimes H^1_{\rm HK})^{N=0,\varphi=p}\ar[r]\ar[d]^{\theta\otimes\iota_{\rm HK}} & 0\\
0\ar[r]& \O/\C_p \ar[r]& \Omega^1 \ar[r] & H^1_{\rm dR}\ar[r] & 0
}$$
\end{theo}
\noindent{\bf Remark}\hskip2mm
We have three proofs of this result.  The proof in~\cite{CDN1} uses the fondamental
exact sequence of sheaves $0\to \Q_p(1)\to(\Bcris^+)^{\varphi=p}\to \widehat\O\to 0$.  
That in~\cite{CDN2}
uses adoc\footnote{An interpolation between ``adic'' and ``ad hoc''.} geometry 
to decompose $Y$ into affinoids with good reduction and open annuli
to give a combinatorial description
of all cohomology groups involved (that is also how Th.\,\ref{cdn7} is proved).
The proof in~\cite{CDN3} uses comparison with syntomic cohomology (as in~\cite{CN,gilles})
and applies to
smooth Stein spaces of arbitrary dimension with a semi-stable model over $\O_{\C_p}$
(this last condition is removed in~\cite{CN4,CN5}).

\vskip3mm
Fix $M$, and take ${\rm JL}(M)$-isotypic components (and factor out the $\check G$-action)
by setting
$$H^1_{\proet}[M]:={\rm Hom}_{\check G}({\rm JL}(M),\varinjlim H^1_{\proet}({\cal M}_{n,\C_p},L(1)))$$
And do the same to define $H^1_{\eet}[M]$, $\O[M]$ and $\Omega^1[M]$

\begin{theo}\phantomsection\label{cdn9}
 {\rm(\cite{CDN1})}
We have a
commutative ${\rm Gal}_{\Q_p}$-equivariant diagram of Fr\'echet spaces:
 $$\xymatrix@R=.3cm@C=.4cm{
0\ar[r] & \O[M]\ar[r]\ar@{=}[d] & H^1_{\proet}[M]\ar[d] \ar[r] &
(\bcris^+\otimes M)^{\varphi=p}\wotimes{\rm LL}(M)^\dual\ar[r]\ar@<-1cm>@{^{(}->}[d] & 0\\
0\ar[r]& \O[M] \ar[r]\ar@{=}[d]& \Omega^1[M] \ar[r] & (\C_p\otimes M_{\rm dR})\wotimes{\rm LL}(M)^\dual\ar[r] & 0\\
0\ar[r]& \O[M] \ar[r]& \C_p\wotimes(\Pi_{M,{\cal L}}^{\rm an})^\dual\ar@{^{(}->}[u] \ar[r] & (\C_p\otimes {\cal L})\wotimes{\rm LL}(M)^\dual\ar[r]\ar@{^{(}->}@<9mm>[u] & 0
}$$
\end{theo}

The top part follows from Th.~\ref{cdn8}, \ref{cdn7} and~\ref{cdn2}.
The bottom part is a consequence of the Breuil-Strauch conjecture proven
by Dospinescu-Le Bras~\cite{DL}.

\vskip3mm
To deduce the first part of Th.\,\ref{cdn4}, one computes that
${\rm Hom}(V_{M,{\cal L}},X)=L$, for $X=\C_p,
(\bcris^+\otimes M)^{\varphi=p}$, from which one infers, by applying
${\rm Hom}(V_{M,{\cal L}},-)$ to the diagram, that
${\rm Hom}(V_{M,{\cal L}},H^1_{\proet}[M])=(\Pi_{M,{\cal L}}^{\rm an})^\dual$.
Finally, Prop.\,\ref{cd} and the fact that $H^1_{\eet}[M]$ is the space of $G$-bounded
vectors in $H^1_{\proet}[M]$ imply 
${\rm Hom}(V_{M,{\cal L}},H^1_{\eet}[M])\hskip.5mm{=}\hskip.5mm\Pi_{M,{\cal L}}^\dual$.

The proof of Th.\,\ref{cdn5} uses the same techniques; that of Th.\,\ref{cdn6}
is more involved as it uses integral $p$-adic Hodge theory (as well
as the description of the rings $R_{M,{\cal B}}$ that can be found in~\cite{CDN6}).

\vskip2mm
\noindent
{\bf Remark}\hskip2mm
We have some very partial results 
for Drinfeld towers associated to other groups than ${\bf GL}_2(\Q_p)$.

(i) In~\cite{CDN3}, we compute the \'etale cohomology of the bottom layer of the tower
for ${\bf GL}_n(F)$, $[F:\Q_p]<\infty$: only Steinberg representations show up on the
${\bf GL}_n(F)$-side and characters on the Galois side.  These results are refined in~\cite{CDN4}
where we use the $\ainf$-cohomology of Bhatt-Morrow-Scholze~\cite{BMS1} to compute the integral
\'etale cohomology. For unexpected applications of these results to algebraic
topology, see~\cite{BSSW2}
(and~\cite{BSSW1}).

(ii) In~\cite{CDN5}, many intermediate results are valid for 
${\bf GL}_2(F)$, $[F:\Q_p]\hskip.5mm{<}\hskip.5mm\infty$.  
These results suggest that the multiplicities
of Galois representations in the \'etale cohomology of the tower associated to
${\bf GL}_n(F)$ are not (dual of) admissible representations of ${\bf GL}_n(F)$
except if $n=2$ and $F=\Q_p$ or at the bottom of the tower.  This creates an interesting
tension with the global situation where multiplicities of Galois representations
in completed cohomology are admissible.


\vskip2mm

\begin{thebibliography}{99}
\bibitem{AW}
K.\,Ardakov, S.\,Wadsley,
On irreducible representations of compact $p$-adic analytic groups,
Ann. of Math.~{\bf 178} (2013), 453--557.


\bibitem{BL94}
{L.\, Barthel}, {R.\,Livn\'e},
Irreducible modular representations of ${\bf GL}_2$ of a local field,
Duke Math. J. {\bf 75} (1994), 261--292.

\bibitem{BSSW1}
T.\,Barthel, T.\,Schlank, N.\,Stapleton, J.\,Weinstein,
On the rationalization of the $K(n)$-local sphere, \url{arXiv:2402.00960},
preprint 2024.

\bibitem{BSSW2}
T.\,Barthel, T.\,Schlank, N.\,Stapleton, J.\,Weinstein,
On Hopkins' Picard group,
\url{arXiv:2407.20958 [math.AT]},
preprint 2024.



\bibitem{BB}
{L.\,Berger}, {C. Breuil}, 
Sur quelques repr\'esentations potentiellement cristallines de $\GL_2(\Q_p)$,
{Ast\'erisque}~{\bf  330} (2010), 155--211.

\bibitem{BMS1} B.\,Bhatt, M.\,Morrow, P.\,Scholze, 
{Integral $p$-adic Hodge Theory},   
Publ. IHES~{\bf 128}  (2018),  219--397.

\bibitem{Br1} {C.\,Breuil}, {Sur quelques repr\'esentations
modulaires et $p$-adiques de ${\GL}_2({\bf Q}_p)$ I},
Comp. Math. {\bf 138} (2003) 165--188.

\bibitem{Br2} {C.\,Breuil}, {Sur quelques repr\'esentations
  modulaires et $p$-adiques de ${\GL}_2(\Q)$ II}, J. Institut
  Math. Jussieu~{\bf 2} (2003), 23--58.

\bibitem{Br3} {C.\,Breuil}, {Invariant\,$\mathcal L$ et s\'erie
sp\'eciale\,$p$-adique},  Ann.\,ENS~{\bf 37} (2004) 559--610.

\bibitem{Br4} {C.\,Breuil},
Introduction g\'en\'erale, Ast\'erisque~{\bf 319} (2008), 1--12.

\bibitem{Br5} {C.\,Breuil},
 Socle\,localement\,analytique\,I,
Ann.\,Inst.\,Fourier\,{\bf 66} (2016), 633--685.

\bibitem{BHS}
C.\,Breuil, E.\,Hellmann, B.\,Schraen, 
A local model for the trianguline variety and applications,
Publ. IHES~{\bf 130} (2019), 299--412.

\bibitem{BHg} {C.\,Breuil}, F.Herzig,
Towards the finite slope part for ${\bf GL}_n$,
Int. Math. Res. Not. IMRN (2020), no.24, 10495--10552.

\bibitem{BD1}
C.\,Breuil, Y.\,Ding, 
Higher ${\cal L}$-invariants for ${\bf GL}_3(\Q_p)$ and local-global compatibility,
Camb. J. Math.~{\bf 8} (2020), 775--951.

\bibitem{BD2}
C.\,Breuil, Y.\,Ding, 
Sur un probl\`eme de compatibilit\'e local-global localement analytique,
Mem. Amer. Math. Soc.~{\bf 290} (2023), no. 1442, vii+152 pp.


\bibitem{6auteurs}
{A.\,Caraiani}, {M.\,Emerton}, {T.\,Gee}, {D.\,Geraghty}, {V.\,\paskunas}, {S.\,W.\,Shin},
Patching and the $p$-adic local Langlands correspondence,
Cambridge J. Math.~{\bf 4} (2016), 197--287.

\bibitem{6auteurs1}
{A.\,Caraiani}, {M.\,Emerton}, {T.\,Gee}, {D.\,Geraghty}, {V.\,\paskunas}, {S.\,W.\,Shin},
 Patching and the  $p$-adic Langlands program for ${\bf GL}_2(\Q_p)$,
Compos. Math.~{\bf 154} (2018), 503--548.


\bibitem{principal}
{P.\,Colmez},
La s\'erie
principale unitaire de ${\bf GL}_2({\bf Q}_p)$,
Ast\'erisque~{\bf 330} (2010), 213--262.

\bibitem{gl2}
{P.\,Colmez},
Repr\'esentations de ${\GL}_2({\bf Q}_p)$ et
$(\varphi,\Gamma)$-modules,
Ast\'erisque~{\bf 330} (2010), 281--509.

\bibitem{poids}
{P.\,Colmez}, Correspondance de Langlands locale $p$-adique et changement de poids,
J.\,EMS~{\bf 21} (2019) 797--838.

\bibitem{adele}
{P.\,Colmez}, 
Exercices ad\'eliques,
\url{arXiv:2402.16231 [math.NT]}, preprint 2024.

\bibitem{CD}
{P.\,Colmez}, {G.\,Dospinescu}, Compl\'etions unitaires de repr\'esentations
de ${\GL}_2(\Q_p)$, Algebra and Number Theory~{\bf 8} (2014), 1447--1519.

\bibitem{CDN1}
{P.\,Colmez}, {G.\,Dospinescu, W.\,Nizio\l},
Cohomologie $p$-adique
de la tour
de Drinfeld, le cas de la dimension~$1$,
J.~AMS~{\bf 33} (2020), 311--362.

\bibitem{CDN2}
{P.\,Colmez}, {G.\,Dospinescu, W.\,Nizio\l},
Cohomologie des courbes analytiques $p$-adiques,
Cambridge J. Math.~{\bf 10} (2022), 511--655.

\bibitem{CDN3}
{P.\,Colmez}, {G.\,Dospinescu, W.\,Nizio\l},
Cohomology of $p$-adic Stein spaces,
Invent. math.~{\bf 219} (2020), 873--985.

\bibitem{CDN4}
{P.\,Colmez}, {G.\,Dospinescu, W.\,Nizio\l},
Integral $p$-adic \'etale Cohomology of Drinfeld symmetric spaces,
Duke Math. J.~{\bf 170} (2021), 575--613.

\bibitem{CDN5}
{P.\,Colmez}, {G.\,Dospinescu, W.\,Nizio\l},
Factorisation de la cohomologie \'etale $p$-adique de la tour de Drinfeld,
Forum Math. Pi~{\bf 11} (2023), e16, 1--62.

\bibitem{CDN6}
{P.\,Colmez}, {G.\,Dospinescu, W.\,Nizio\l},
Correspondance de Langlands locale $p$-adique et anneaux de Kisin,
Acta Arithmetica~{\bf 208} (2023), 101--126.


\bibitem{CDP}
{P.\,Colmez}, {G.\,Dospinescu}, {V.\,\paskunas},
The $p$-adic local Langlands correspondence for ${\GL}_2(\Q_p)$,
Cambridge J. Math.~{\bf 2} (2014), 1--47.

\bibitem{CF} 
{P.\,Colmez}, {J.-M.\,Fontaine},
Construction des repr\'esentations $p$-adiques semi-stables,
Invent. math.~{\bf 140} (2000), 1--43.

\bibitem{CN} {P.\,Colmez}, {W.\,Nizio\l},
{Syntomic complexes and $p$-adic nearby cycles},
Invent. math.~{\bf 208} (2017) 1--108.

\bibitem{CN4}
{P.\,Colmez}, {W.\,Nizio\l},
On the cohomology of $p$-adic analytic spaces, I:  The basic comparison theorem,
J. Algebraic Geometry~{\bf 34} (2025), 1--108.

\bibitem{CN5}
{P.\,Colmez}, {W.\,Nizio\l},
On the cohomology of $p$-adic analytic spaces, II:  The $C_{\rm st}$-conjecture,
\url{arXiv:2108.12785 [math.AG]}, preprint 2021, to appear in Duke Math. J. (2025).

\bibitem{CW}
{P.\,Colmez}, {S.\,Wang},
Une factorisation de la cohomologie compl\'et\'ee et du syst\`eme de Beilinson-Kato,
\url{arXiv:2104.09200 [math.NT]}, preprint 2021, nouvelle version 2024.


\bibitem{Dosp1}
{G.\,Dospinescu}, Actions infinit\'esimales dans la correspondance
de Langlands locale $p$-adique pour ${\GL}_2(\Q_p)$,
Math.~Ann.~{\bf 354} (2012), 627--657.

\bibitem{DL}
{G.\,Dospinescu}, {A.-C.\,Le Bras},
Rev\^etements du demi-plan de Drinfeld et correspondance de
Langlands locale $p$-adique, Ann.~of Math.~{\bf 186} (2017), 321--411.

\bibitem{DS}
{G.\,Dospinescu}, {B.\,Schraen},
Endomorphism algebras of $p$-adic representations of $p$-adic Lie groups, Representation Theory~{\bf 17}
 (2013), 237--246.

\bibitem{Drinfeld}
{V.\,Drinfeld},
 Elliptic modules, Math. Sb.~{\bf 94} (1974), 594-627.

 \bibitem{Drinfeld1}
{V.\,Drinfeld}, Coverings of $p$-adic symmetric regions, Funktsional. Anal. i Prilozhen., (1976), 29-40; Funct. Anal. Appl., (1976), 107-115.


\bibitem{Em05}
{M.\,Emerton},
$p$-adic $L$-functions and unitary completions of representations of $p$-adic reductive groups,
 Duke Math. J.~{\bf 130} (2005), 353--392.

\bibitem{Em06}
{M.\,Emerton},
On the interpolation of systems of eigenvalues attached to automorphic Hecke eigenforms,
Invent. math.~{\bf 164} (2006), 1--84.

\bibitem{Em06b} {M.\,Emerton},
A local-global compatibility conjecture in the $p$-adic Langlands programme
for ${\GL}_{2/{\Bbb Q}}$,
Pure Appl. Math. Q.~{\bf 2} (2006), 279--393.

\bibitem{Em08}
{M.\,Emerton}, Local-global compatibility in the $p$-adic
Langlands programme for ${\GL}_{2,\Q}$, preprint 2008!

\bibitem{Em14}
{M.\,Emerton},
Completed cohomology and the p-adic Langlands program,
{\it Proceedings of the 2014 ICM, Seoul}, \url{http://www.icm2014.org/en/vod/proceedings.html}.

\bibitem{EGH}
{ M.\,Emerton, T.\.Gee, E.\,Hellmann},
An introduction to the categorical $p$-adic Langlands program,
\url{arXiv:2210.01404 [math.NT]}, preprint 2022.


\bibitem{Fo91}  {J.-M.\,Fontaine},
Repr\'esentations $p$-adiques des
corps locaux, in  ``{\em The Grothendieck Festschrift}'',
vol 2, Prog. in Math.~{\bf 87}, 249--309, Birkh\"auser 1991.

\bibitem{gilles} S.\,Gilles, 
{Morphismes de p\'eriodes et cohomologie syntomique}.      
 Algebra Number Theory~{\bf 17} (2023), 603--666.

\bibitem{Ki-FM}
 {M.\,Kisin},
The\,Fontaine-Mazur\,conjecture\,for\,${\bf GL}_2$,
 J.\,AMS~{\bf 22}\,(2009),\,641--690.

\bibitem{pan1}
{L.\,Pan},
The Fontaine-Mazur conjecture in the residually reducible case,
J. AMS~{\bf 35} (2022), 1031--1169.

\bibitem{pan2}
{L.\,Pan},
On locally analytic vectors of the completed cohomology of modular curves,
 Forum Math. Pi~{\bf 10} (2022), e7, 82 pp.


\bibitem{Pas}
{V.\,\paskunas},
The image of Colmez's Montreal functor, Publ. IHES~{\bf 118} (2013), 1--191.

\bibitem{Pas1}
{V.\,\paskunas},  
{Blocks for mod $p$ representations of ${\bf GL}_2(\Q_p)$}, 
{\it Automorphic forms and Galois representations}, Vol. 2, 231--247,
London Math. Soc. Lecture Note Ser., 415, Cambridge Univ. Press, 2014.

\bibitem{PT}
{\sc V.\,\paskunas, S.-N.\,Tung},
Finiteness properties of the category of mod $p$ representations of ${\rm GL}_2(\Q_p)$.
Forum Math. Sigma~{\bf 9} (2021), e80, 39 pp.

\bibitem{ST1} {P.\,Schneider}, {J.\,Teitelbaum},
Locally analytic
  distributions and $p$-adic representation theory, with applications
  to ${\GL}_{2}$, J. AMS~{\bf15} (2002), 443--468.


\bibitem{ST2} {P.\,Schneider}, {J.\,Teitelbaum},
(with an appendix by
{ D. Prasad}), {$U(\mathfrak g)$-finite locally analytic
    representations}, Representation Theory~{\bf 5} (2001), 111--128.

\bibitem{ST3} {P.\,Schneider}, {J.\,Teitelbaum},
{Banach space
  representations and Iwasawa theory}, Israel J. Math.~{\bf 127} (2002),
  359--380.

\bibitem{ST4} {P.\,Schneider}, {J.\,Teitelbaum},
Algebras of $p$-adic distributions and admissible representations,
 Invent. Math.~{\bf  153} (2003), 145--196.


\bibitem{VanH}
{A.\,Vanhaecke},
Cohomologie de syst\`emes locaux $p$-adiques sur les rev\^etements du demi-plan de Drinfeld,
\url{arXiv:2405.10048 [math.NT]}, preprint 2024.



\end{thebibliography}
\end{document}